\documentclass[letterpaper,12pt,leqno]{article}

\usepackage{fullpage}

\usepackage{hyperref}

\usepackage{amsmath}
\usepackage{amsfonts}
\usepackage{amssymb}
\usepackage{tabularx}
\usepackage{color}
\usepackage{fancyhdr}
\usepackage[all]{xy}
\usepackage{enumerate}

\def\theoremparent{section}
\usepackage{notes}

\numberwithin{equation}{section}
\numberwithin{theorem}{subsection}

\def\b{\mathrm{b}}
\def\bb{\mathfrak{b}}

\def\Br{\mathrm{Br}}

\def\cD{\mathcal{D}}
\def\cE{\mathcal{E}}

\def\d{\partial}

\def\Der{\mathrm{Der}}

\def\E{\mathcal E}

\def\G{\mathcal{G}}

\def\In{\mathrm{In}}

\def\MPP{\mathrm{MPP}}

\def\op{\mathrm{op}}

\def\Out{\mathrm{Out}}

\def\PProp{\mathrm{PProp}}

\def\proof{\noindent{\em Proof:}\ }

\def\qed{\hfill\lower 1em\hbox{$\square$}\vskip 1em}
\def\Rees{\mathrm{Rees}}

\def\toto{\text{\ \raise0.2em\hbox to 0pt {$\to$}\lower0.2em\hbox{$\to$}\ }}

\def\Vect{\mathrm{Vect}}

\def\lsqcup#1{{\mathop{\ \sqcup\ }\limits^{\scriptscriptstyle <}_{#1}}}
\def\sh{\ \hat{*}\ }

\begin{document}

\title{Planar Prop of Differential Operators}
\author{Slava Pimenov}
\date{\today}
\titlepage
\maketitle

\tableofcontents
\vskip 5em

\setcounter{section}{-1}
\section{Introduction}
There are many different approaches to generalize the notion of differential calculus to the noncommutative setting, extending various
aspects of the commutative theory (such as \cite{Connes}, \cite{CQ2}, \cite{LR}, \cite{Sa}, \cite{TT} to name a few). They
produce distinct albeit interconnected and equally meaningful theories. In this paper we propose another approach to extend the
notion of differential operators to the setting of associative algebras.

We believe that a ``correct'' notion of the differential operators should satisfy the following principles and we
will be using them as guidelines for our definition:
\begin{enumerate}[a)]
\item the construction of differential operators should respect quasi-isomorphisms between dg-algebras,
\item there should be a natural symbol map from differential operators to polyvector fields,
\item the role of the associative algebra structure should be made explicit in the construction.
\end{enumerate}

Since the notion of differential operators is not a functor, the first principle is better applied to the corresponding notion
of $D$-modules. Basically it says that the definition should be readily extendible to dg-algebras and it should be a Quillen
functor, so that its derived functor can be computed using cofibrant resolutions of associative dg-algebras. In a sense,
this can be thought of as a dg-version of Kashiwara equivalence theorem.

Principle (b) essentially means that whatever we define has a right to be called differential operators. For our construction
in this paper we convert it to a precise statement (proposition \ref{symbol.map}).

The reason for the third principle is the fact that differential operators are closely related to the deformation theory.
Instead of associative algebras one could consider the category of $\cP$-algebras with the algebraic structure given for instance
by some operad $\cP$ and ask what is the corresponding notion of $\cP$-differential operators and its connection to the
deformation theory of $\cP$-algebras. Of course one expects that
when applied to the operad $\mathcal Com$ of commutative algebras one recovers the classical notion of differential operators.
In this paper we will not be using this principle explicitly in the definition, however in section \ref{sect-assoc-families}
we construct another planar prop $\cE(A)$ from automoprhisms of the trivial formal associative deformation of $A$
with a natural map to the planar prop of multi-differential operators $\cD(A)$.

\begin{nparagraph}
As the basis for our definition we use the following identities satisfied by the differential operators of a commutative
algebra $A$. For an operator $P$ of order $1$
we have
$$
P(ab) = P(a)b + aP(b) - aP(1)b,
$$
for an order $2$ operator we have
$$
P(abc) = P(ab)c + aP(bc) + bP(ac) - P(a)bc - aP(b)c - abP(c) + abcP(1),
$$
and so forth. While the first relation can be applied without change to an associative algebra, the relations for higher order
operators can not. In order to address this issue our operator $P\from A \to A$ comes with ``companions'' $P_k\from A^{\tensor k} \to A^{\tensor k}$,
which are also poly-differential maps. Each companion measures the defect of its predecessor to be a poly-derivation.

The spaces $D(A)$ (proposition \ref{assoc.DA}) and $\cD(A)$ (theorem \ref{thm.planar.prop}) that we define are analogs of the Rees algebra
$\Rees_F(D)$ associated to the standard filtration on differential operators in commutative setting.

\end{nparagraph}

\begin{nparagraph}[Outline of the paper.]
In the first section we define an algebraic structure that we call the planar prop. It is part of the menagerie of similar
structures, such as operads, properads, half-props, props, etc. In section \ref{sect-graph-def} we give the precise definition
using language of graph-based Feynman categories to set up the notation and clarify the details. Then in section \ref{sect-elem-def}
we convert this definition of a planar prop into a more elementary one, as a collection of spaces equipped with two operations
satisfying certain explicit relations. This alternative definition will be used to show that the space of multi-differential
operators has a structure of a planar prop.

We begin section \ref{sect-diff-def} by introducing some notations and constructions concerning the category $\Delta$ of totally
ordered finite sets that will be carrying the combinatorics of the definition of differential operators. Then we proceed with
giving the main definition of operators at $f$ for a map of associative algebras $f\from A \to B$, which can be thought as
operators on $A$ with values in $B$. Finally, we extend this to operators with multiple inputs and outputs and call the resulting
space the multi-differential operators $\cD(A)$.

Section \ref{sect-planar-prop} is dedicated to showing that the constructed space $\cD(A)$ equipped with bigrading by the number of
inputs and outputs forms a planar prop. We also observe that if we assign a genus to an operator, which is basically the
difference between its order and the number of outputs, we in fact obtain the modular planar prop structure.

In section \ref{sect-assoc-families} we consider trivial associative families of algebras with fiber $A$ over completions
of free associative algebras. We use automorphisms of such families as a source of differential operators. In fact we
construct another planar prop $\cE(A)$ with a map to $\cD(A)$. The reason for this construction is to follow through with
the third guiding principle. This is also the reason to consider operators with multiple inputs and outputs.

Finally, the section \ref{sect-symbol} is dedicated to the case of formally smooth associative algebras. We show that in this
case there is a surjective symbol map from $\cD(A)$ to poly-derivations and we identify its kernel as the image of simplicial degeneracy
maps. We also show that for formally smooth algebras the image of the natural map $\cE(A) \to \cD(A)$ is formed by operators
of totally positive genus.
\end{nparagraph}

The author would like to thank Mikhail Kapranov for his useful remarks and suggestions in preparation of this paper.

\vfill\eject

\section{Planar Prop}

We define the notion of a planar prop analogously to that of a prop (see for example \cite{Vallette}) except we do not keep
track of the action of symmetric groups on the inputs and outputs of operations, and we only allow compositions along
planar graphs. First, we give the definition using formalism of Feynman categories (\cite{Kaufmann}).

\subsection{Graph-based definition}
\label{sect-graph-def}

\begin{nparagraph}
A {\em directed} graph $G$ consists of finite sets of vertices $V$ and half-edges $E$, equipped with
\begin{itemize}
\item a map $\alpha\from E \to \{\In, \Out\}$, decorating half-edges as either inputs or outputs,
\item the incidence map $i\from E \to V$ describing how the half-edges are attached to vertices,
\item the involution $\epsilon\from E \to E$ describing how the half-edges are glued together to form edges. This involution
is compatible with the decoration $\alpha$, in the sense that outputs are glued with inputs.
\end{itemize}

For each vertex $v \in V$ we denote $\In(v)$ and $\Out(v)$ the sets of incident half-edges marked as inputs and outputs
respectively. A half-edge fixed by involution $\epsilon$ and marked as input (output) is called an input (respectively output) of $G$.
The rest of the $\Z/2$-orbits with the free action of the involution are called internal edges. In other words an internal edge
is an ordered pair $(e_1, e_2)$ such that $e_1 \in \Out(v_1)$ and $e_2 \in \In(v_2)$ for some vertices $v_1$ and $v_2$.
We say that this edge goes from $v_1$ to $v_2$.

A directed cycle in $G$ is a sequence of internal edges $(e_k, f_k)$, $1 \le k \le n$, such that $i(f_k) = i(e_{k+1})$
for $1 \le k \le (n-1)$ and $i(f_n) = i(e_1)$.
$G$ is called {\em acyclic} if it has no directed cycles.

We say that $G$ is a {\em ribbon} graph if the sets $i^{-1}(v)$ are cyclically ordered for every $v \in V$.
Furthermore, we say that $G$ is {\em planar} if
\begin{itemize}
\item the cyclic order on $i^{-1}(v)$ is such that all inputs (and therefore outputs) are gathered together.
In other words, if $v$ has both inputs and outputs there are $i \in \In(v)$ and $o \in \Out(v)$, such that
\begin{equation}
\label{cyc.order}
\forall x \in \In(v): i \preceq x \prec o, \quad\text{and}\quad \forall y \in \Out(v): o \preceq y \prec i.
\end{equation}
\item the geometric realization of $G$ can be embedded into a unit disk in $\R^2$ respecting the cyclic order at each vertex,
with the free ends of inputs and outputs of $G$ lying on the boundary of the disk, in a way that the induced cyclic order
satisfies \ref{cyc.order}
\end{itemize}

We say that graph $G$ is {\em essential} if every vertex has at least one input and at least one output.

All graphs are assumed to be directed acyclic planar graphs unless stated otherwise. We denote the set of such graphs by $\G$.

\end{nparagraph}

\begin{remark}
Our notion of a planar graph differs from the usual notion by imposing condition \ref{cyc.order}. If $G$ is essential, it is equivalent to
providing linear orders on the (non-empty) sets $\In(v)$ and $\Out(v)$. Alternatively, one could think of such a graph as having a
global direction, in the sense that there exists an embedding as above, with additional property that all edges go in
a prescribed direction, say downward (see lemma \ref{level.graph}).

\end{remark}

\begin{nparagraph}[Genus markings.]
\label{genus.markings}
Let $G$ be a graph as above, and consider a function $\gamma\from V \to \Z$ that we will refer to as genus marking of
vertices. We allow negative integers here, since our graphs may not be connected. We define the genus of $G$ as
$$
\gamma(G) = \dim H_1(G) - \dim H_0(G) + 1 + \sum_{v \in V} \gamma(v).
$$
One can check that this definition is compatible with the composition of graphs. Fix vertex $v \in V$ and take a genus marked graph $H$, such
that number of inputs and outputs of $H$ are the same as inputs and outputs of $v$, and $\gamma(H) = \gamma(v)$.
The composition $H \circ_v G$ is the graph obtained from $G$ by replacing vertex $v$ with $H$ and connecting
inputs and outputs of $v$ to inputs and outputs of $H$ preserving the cyclic ordering. Then $\gamma(H \circ_v G) = \gamma(G)$.

We will write $\G_m$ for the set of genus marked graphs.
\end{nparagraph}

\begin{nparagraph}[Feynman category.]
A Feynman category is a triple $(\mathcal V, \mathcal F, i\from \mathcal V \into \mathcal F)$, where $\mathcal V$ is a groupoid,
$\mathcal F$ is a symmetric monoidal category and $i$ is an inclusion functor, that satisfy the following conditions.
\begin{enumerate}[a)]
\item Isomorphism condition: the induced functor $i^{\tensor} \from \mathcal V^{\tensor} \to \Iso(\mathcal F)$ is an equivalence.
Here $\mathcal V^{\tensor}$ denotes the symmetric monoidal category generated by $\mathcal V$.
\item Hereditary condition: the induced functor $(\Iso(\mathcal F \downarrow \mathcal V))^{\tensor} \to \Iso(\mathcal F \downarrow \mathcal F)$
is an equivalence.
\item Size condition: for any $v \in \mathcal V$, $(\mathcal F \downarrow v)$ is a small category.
\end{enumerate}

As is customary with graph-based Feynman categories we put $\mathcal V = \mathcal Cor$, the groupoid of corollas,
i.e. graphs with a single vertex. Isomorphism class of an object $v \in \mathcal Cor$ is determined by a pairs of non-negative
integers $(m, n)$, where $m = |\Out(v)|$ and $n = |\In(v)|$. We set $\Aut(v_{m,n}) = 1$ if $mn \neq 0$ and
$\Aut(v_{n, 0}) = \Aut(v_{0,n}) = \Z/n$.

Objects of $\mathcal F$ are collections
of corollas (i.e graphs without internal edges), the monoidal structure is given by disjoint union. In virtue of the
hereditary condition above, a map in $\mathcal F$ to a collection of corollas is completely determined by a map to each
of the corolla in the collection. We set $\Hom_{\mathcal F}(X, v)$ to be the set of graphs $G$, such that number of inputs
and outputs of $G$ are the same as for corolla $v$, and vertices of $G$ together with their inputs and outputs are in
bijection with the corollas in $X$. The composition of morphisms is given by the composition of graphs.

Similarly, we form the genus marked Feynman category $\mathcal F_m$ by taking $\mathcal V = \mathcal Cor \times \Z$ the
groupoid of corollas marked by an integer, and $\Hom_{\mathcal F_m}(X, v)$ to be the set of suitable genus marked graphs.

\end{nparagraph}

\begin{definition}
A planar prop is a symmetric monoidal functor $\cP\from \mathcal F \to \Vect$. A modular planar prop is a symmetric monoidal
functor $\mathcal Q\from \mathcal F_m \to \Vect$.
\end{definition}

The term {\em modular} is used to indicate that this notion is basically a planar analog of the modular operad of
Getzler-Kapranov (\cite{Getz.Kap}).

Since a corolla $v$ is completely determined by a pair of non-negative integers $(m, n)$, where $m = |\Out(v)|$, $n = |\In(v)|$, the
restriction of a planar prop $\cP$ to the groupoid of corollas gives a collection of vector spaces $\cP(m, n)$. Informally, the
planar prop structure prescribes how to compose elements of $\cP(m, n)$ along planar graphs.

\begin{nparagraph}
We denote the category of planar props by $\PProp$ and the category of modular planar props by $\MPP$. Consider the monoidal
structure $\tensor$ on planar props given by the coproduct in $\PProp$. The unit for this monoidal structure is given by the planar
prop $\mathcal I$, with
$$
\mathcal I(m, n) = \begin{cases}
k,\quad\text{if $m = n$},\\
0,\quad\text{otherwise}.
\end{cases}
$$
We say that a planar prop $\cP$ is {\em unital} if it admits a unit map $u\from \mathcal I \to \cP$, i.e. the composition
$$
\xymatrix@C=3em{
\cP \ar^-\isom[r] & \mathcal I \tensor \cP \ar^-{u \tensor \id}[r] & \cP \tensor \cP \ar[r] & \cP
}
$$
is identity.

\end{nparagraph}

\begin{nparagraph}
From the definition of the unital planar prop $\cP$ we immediately see that
\begin{enumerate}[a)]
\item $\cP(0, 0)$ is a commutative algebra,
\item $\cP(m, n)$ are $\cP(0, 0)$-modules for all $m, n \in \Z_{\ge 0}$,
\item $\cP(1, 1)$ is an associative $\cP(0, 0)$-algebra, and $\cP(0, 0)$ is contained in its center,
\item $\cP(n, 0)$ and $\cP(0, n)$ are equipped with the action of the cyclic group $\Z/n$, and an invariant pairing
$\cP(0, n) \tensor \cP(n, 0) \to \cP(0, 0)$.
\end{enumerate}

A unital planar prop $\cP$ is called {\em reduced} if $\cP(n, 0) = \cP(0, n) = 0$ for $n \ge 1$ and $\cP(0, 0) = k$.
In the rest of the paper we will be mostly concerned with reduced planar props.
\end{nparagraph}

\vskip 5em
\subsection{Elementary definition}
\label{sect-elem-def}
Let us give an alternative definition of a planar prop as a collection of vector spaces and operations.

\begin{definition}
\label{def.elem.prop}
An elementary prop $\cP$ is a collection of vector spaces $\cP(m, n)$, $m, n \in \Z_{\ge 1}$, and $\cP(0, 0) = k$,
with a vertical unit element $u \in \cP(1, 1)$, horizontal unit $u^0 \in \cP(0, 0)$, and two operations: horizontal composition
$$
\circ_h \from \cP(m_1, n_1) \tensor \cP(m_2, n_2) \to \cP(m_1 + m_2, n_1 + n_2),
$$
and vertical composition
$$
\circ_v \from \cP(m, n) \tensor \cP(n, k) \to \cP(m, k).
$$
Operations $\circ_h$ and $\circ_v$ are associative and satisfy the following relations. Let $a \in \cP(m, n)$
and $b \in \cP(p, q)$, and denote by $u^k \in \cP(k, k)$ the horizontal composition of $k$ copies of $u$, for $k \ge 1$.
\begin{enumerate}[a)]
\item Horizontal unit:
$$
u^0 \circ_h a = a = a \circ_h u^0.
$$
\item Compatibility:
$$
(a \circ_h u^p) \circ_v (u^n \circ_h b) = a \circ_h b = (u^m \circ_h b) \circ_v (a \circ_h u^q).
$$
\end{enumerate}

\end{definition}

A morphism between two elementary props is a map between corresponding collections preserving units and the two
compositions.

\begin{definition}
The expression of the type
$$
(u^{i_1} \circ_h a_1 \circ_h u^{j_1}) \circ_v \cdots \circ_v (u^{i_r} \circ_h a_r \circ_h u^{j_r})
$$
with $a_i \in \cP(m_i, n_i)$
is said to be in the normal form if for any $1 \le k \le r$ and $1 \le l \le (k-1)$
$$
i_k < \min(i_l, \ldots i_{k-1}) \quad \text{implies} \quad i_k + m_k > \min(i_l, \ldots i_{k-1}).
$$
\end{definition}

\begin{proposition}
\label{elem.prop}
The category of elementary props is equivalent to the category of reduced unital planar props.
\end{proposition}

Before proving the proposition let us show that our planar graphs can be embedded into $\R^2$ in a certain way.

\begin{lemma}
\label{level.graph}
Every essential graph $G$, admits a level embedding into $\R^2$, i.e. an embedding such that
\begin{itemize}
\item all vertices of $G$ lie in $\R \times \Z \subset \R^2$,
\item every line $\R \times \{n\}$ contains at most one vertex,
\item every edge goes downwards.
\item there is $N \in \Z$, such that all vertices lie in $\R \times (-N, N)$, free ends of inputs of $G$ lie on
$\R \times \{N\}$ and free ends of outputs lie on $\R \times \{-N\}$.
\end{itemize}
Moreover, there is a canonical linear ordering on the set of vertices of $G$, and this embedding can be constructed in
such a way that the ordering of vertices by their levels coincides with this canonical order.
\end{lemma}

\proof
We prove the statement by induction on the maximal number of vertices in a connected component of $G$. Clearly, if each
connected component of $G$ admits the required embedding then so does the entire $G$.
If $G$ has only one vertex there is nothing to show.
Otherwise, since $G$ has no directed cycles, there exists a vertex $v$ such that all outputs of $v$ are outputs of $G$.
Take the first such vertex with respect to the linear order on $\Out(G)$, and construct a new graph $G'$ first by removing
vertex $v$ from $G$ along with all half-edges in $\Out(v)$ as well as half-edges in $\In(v)$ which are part of an internal edge
of $G$. Then for all half-edges $e \in \In(v) \cap \In(G)$ we add a new vertex $v_e$ and a half-edge $e'$, so that
$\In(v_e) = \{e\}$ and $\Out(v_e) = \{e'\}$.

By construction, $G'$ is again an essential graph and since each new vertex $v_e$ belongs to its own connected component the maximal
number of vertices of a connected component of $G'$ is strictly less than that of $G$. Therefore, it admits a level embedding into $\R \times [-N', N']$. Let us extend
this embedding to $G$: we place $v$ on $\R \times \{-N'\}$, connect its inputs to appropriate outputs of $G'$ by line segments
going down from $\R \times \{-N' + 1\}$ to $v \in \R \times \{-N'\}$ (this can be done because of the cyclic order condition
\ref{cyc.order}), then embed outputs of $v$ as line segments going
from $v$ to $\R \times \{-N' - 1\}$, and finally extend the remaining outputs of $G'$ down to $\R \times \{-N' - 1\}$.
This gives us a level embedding of $G$ into $\R \times [N'+1, -N'-1]$.

\qed

\begin{nparagraph}[Proof of the proposition.]
In one direction the statement is straightforward. Let $\cP$ be a reduced unital planar prop, put $u$ to be the image of
$1 \in \mathcal I(1, 1)$ and $u^0$ the image of $1 \in \mathcal I(0, 0)$ under the unit map $\mathcal I \to \cP$. Define
horizontal composition $a \circ_h b$ by the graph
$$
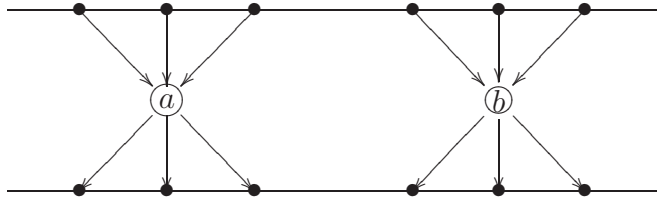
\xymatrix{
\ar@{-}[0,8] & *=0{\bullet} \ar[dr] & *=0{\bullet} \ar[d] & *=0{\bullet} \ar[dl] & & *=0{\bullet} \ar[dr] & *=0{\bullet} \ar[d] & *=0{\bullet} \ar[dl] & \\
&& *+[Fo]{a} \ar[dl] \ar[d] \ar[dr] && && *+[Fo]{b} \ar[dl] \ar[d] \ar[dr] && \\
\ar@{-}[0,8] & *=0{\bullet} & *=0{\bullet} & *=0{\bullet} & & *=0{\bullet} & *=0{\bullet} & *=0{\bullet} &
}
$$
and the vertical composition $a \circ_v b$ by
$$
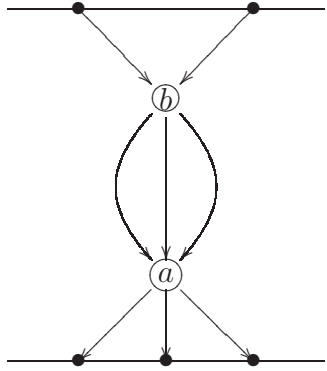
\xymatrix{
\ar@{-}[0,4] & *=0{\bullet} \ar[dr] & & *=0{\bullet} \ar[dl] & \\ 
&& *+[Fo]{b} \ar@(dl,ul)[dd] \ar[dd] \ar@(dr,ur)[dd] && \\
&&&& \\
&& *+[Fo]{a} \ar[dl] \ar[d] \ar[dr] && \\
\ar@{-}[0,4] & *=0{\bullet} & *=0{\bullet} & *=0{\bullet} &
}
$$
The relations for $\circ_v$ and $\circ_h$ immediately follow from the properties of the unit map $\mathcal I \to \cP$ and
the fact the corresponding compositions are represented by identical graphs.

Now, let $\cP$ be an elementary prop, we want to show that it has a structure of a reduced unital planar prop.
Since the prop is reduced it is enough to define operations corresponding to essential graphs. According to lemma \ref{level.graph}
the composition along such a graph can be written in the normal form uniquely determined by the graph. Let $G$ be an essential
graph, and $(v_1, \ldots v_r)$ be the set of its vertices ordered as in the lemma. Denote $m_i = |\Out(v_i)|$ and $n_i = |\In(v_i)|$,
and let $a_i \in \cP(m_i, n_i)$. Then
$$
\cP(G)(a_1, \ldots a_r) = (u^{i_1} \circ_h a_1 \circ_h u^{j_1}) \circ_v \cdots \circ_v (u^{i_r} \circ_h a_r \circ_h u^{j_r})
$$
is in the normal form.

It remains to show that this assignment is well defined, in other words that any expression involving $\circ_h$ and $\circ_v$
can be reduced to the normal form using relations in \ref{def.elem.prop}. First, using associativity of compositions and the
compatibility relation we rewrite it as a vertical composition of multiple terms of the form $(u^i \circ_h a \circ_h u^j)$.
If this expression is not in the normal form, then again using compatibility relation we can push the offending term
towards the beginning of the product. Repeating this process we reduce the expression to the normal form.

Finally, observe that the horizontal and vertical units determine a map $\mathcal I \to \cP$, and combining relations (a)
and (b) we find that it is a unit map.

\qed

\end{nparagraph}

\begin{example}[Braid prop.]
As an illustration of the notion we give the following example. Denote $\Br_n$ the braid group on $n$ strands, and put
$$
\mathcal Br(m, n) = \begin{cases}
k[\Br_n],\quad\text{if $m = n$}\\
0,\quad\text{otherwise}.
\end{cases}
$$
This is a unital planar prop generated by a single element $s \in \mathcal Br(2, 2)$, modulo relation
$$
(s \circ_h u) \circ_v (u \circ_h s) \circ_v (s \circ_h u) = (u \circ_h s) \circ_v (s \circ_h u) \circ_v (u \circ_h s).
$$

\end{example}

\vfill\eject

\section{Differential operators}

\subsection{Definition and basic properties}
\label{sect-diff-def}

\begin{definition}
An ordered partition of $n \in \Z_{\ge 0}$ of size $d$ is a sequence $\lambda_\bullet = (\lambda_1, \ldots \lambda_d)$,
such that each $\lambda_i \ge 0$ and $\sum_{i=1}^{d} \lambda_i = n$.

A partition $\lambda_\bullet$ is said to be degenerate if either $n > 0$ and at least one of $\lambda_i = 0$, or $n = 0$ and $d > 1$.
\end{definition}

One can think of such a partition as a decomposition of the integral interval $I = [1,n] \cap \Z$ into a disjoint union of
intervals (possibly empty) $I = \bigcup_{i=1}^d I_i$, where $|I_i| = \lambda_i$. We will write $|\lambda_\bullet|$
for the size of the partition.
Denote by $\Lambda_d(n)$ the set of all ordered partitions of $n$ of size $d$, and $\Lambda_d^\circ(n)$ the subset of non-degenerate partitions.

Let $\lambda_\bullet$ be a partition of $n$ of size $k$ and $\mu_\bullet$ a partition of some $\lambda_i$ of size $l$.
Then we define the composition
$$
(\mu \circ_i \lambda)_\bullet = (\lambda_1, \ldots \lambda_{i-1}, \mu_1, \ldots \mu_l, \lambda_{i+1}, \ldots \lambda_k).
$$

We say that a partition $\lambda'$ is a {\em refinement} of $\lambda$ if $\lambda'$ can be obtained from $\lambda$ by a
sequence of compositions
$$
\lambda' = \mu^{(1)} \circ_{i_1} ( \mu^{(2)} \circ_{i_2} \cdots (\mu^{(p)} \circ_{i_p} \lambda) \cdots ).
$$
In other words, every interval in partition $\lambda'$ is contained in some interval in $\lambda$.
A refinement is said to be {\em non-degenerate} if all $\mu^{(j)}$ are non-degenerate.

Furthermore, we form the category of non-degenerate partitions $\ul\Lambda^\circ(n)$ by setting the set of morphisms
$\Hom_{\ul\Lambda^\circ}(\lambda', \lambda)$ consist of a single element if $\lambda'$ is a refinement of $\lambda$ and
an empty set otherwise.

\begin{remark}
\label{rem.Delta}
Let $\Delta$ be the indexing category for augmented simplicial sets, i.e. the category of totally ordered finite sets and
order preserving maps between them. We denote by $[n]$ the set with $n$ elements (the empty set if $n = 0$).
It is clear that we have $\Lambda_d(n) = \Hom_\Delta([n], [d])$ and
$\Lambda_d^\circ(n) = \mathrm{Epi}_\Delta([n], [d])$, the set of surjective maps from $[n]$ to $[d]$ in $\Delta$.
Consider two partitions $\lambda \in \Lambda_d(n)$ and $\lambda' \in \Lambda_{d'}(n)$ and corresponding maps $[n] \to [d]$ and
$[n] \to [d']$ in $\Delta$. It is clear that $\lambda'$ is a refinement of $\lambda$ if and only if there exists a surjective
map $\rho\from [d'] \epi [d]$ forming the commutative triangle
$$
\xymatrix{
[n] \ar_{\lambda}[dr] \ar^{\lambda'}[rr] && [d'] \ar@{->>}^{\rho}[dl] \\
& [d]. &
}
$$
If it is a non-degenerate refinement then such $\rho$ is unique. In fact the uniqueness of $\rho$ can be used to characterize
non-degenerate refinements: $\lambda'$ is a non-degenerate refinement of $\lambda$ if and only if for any partitions $\mu$, $\nu$
there exists unique map $\rho$ associated to the refinement $(\mu, \lambda', \nu)$ of concatenation $(\mu, \lambda, \nu)$
(see the next paragraph for the definition of concatenation).

So the category $\ul\Lambda^\circ(n)$
is equivalent to the full subcategory of $([n] \downarrow \Delta)$ with objects being surjective maps $[n] \epi [d]$ for some
$[d] \in \Delta$.

Similarly, we define category $\ul\Lambda(n)$ with objects $\Lambda_d(n)$ for all $d$ and morphisms
$$
\Hom_{\ul\Lambda(n)}(\lambda', \lambda) = \mathrm{Epi}_{[n] \downarrow \Delta}(\lambda', \lambda).
$$

\end{remark}

\begin{nparagraph}
Let us introduce some additional terminology regarding ordered partitions that will be needed in the future.
For two partitions $\lambda \in \Lambda_p(n)$ and $\mu \in \Lambda_q(m)$ we define the concatenation
$$
(\lambda, \mu) = (\lambda_1, \ldots, \lambda_p, \mu_1 \ldots, \mu_q) \in \Lambda_{p+q}(m+n).
$$
Consider a refinement $\nu$ of the concatenation $(\lambda, \mu)$, we will denote by $\nu|_\lambda$ and $\nu|_\mu$
the refinements of $\lambda$ and $\mu$ respectively obtained by restriction of $\nu$ to corresponding subpartition.

Let us express this using vocabulary of the category $\Delta$ introduced in remark (\ref{rem.Delta}). For two objects
$[m], [n] \in \Delta$ we define the {\em ordered coproduct} $[m] \lsqcup{} [n]$ as the initial object in the category
of pairs $(\alpha\from [m] \to [p], \beta\from [n] \to [p])$, such that for any $x \in [m]$, and $y \in [n]$ we have
$\alpha(x) < \beta(y)$ in $[p]$. In other words $[m] \lsqcup{} [n]$ is isomorphic to $[m + n]$ with $[m]$ embedded as the
first $m$ elements and $[n]$ as the last.

We also have the relative version of this: let $f\from [m] \to [q]$ and $g\from [n] \to [q]$, we define $[m] \lsqcup{[q]} [n]$
as the initial object in the category of pairs of maps $(\alpha\from [m] \to [p], \beta\from [n] \to [p])$ in $\Delta \downarrow [q]$,
such that for any $x \in [m]$ and $y \in [n]$ with $f(x) = g(y)$ we have $\alpha(x) < \beta(y)$ in $[p]$.

Now, the concatenation of two partitions $\lambda \in \Lambda_p(n)$ and $\mu \in \Lambda_q(m)$ is represented by the ordered
coproduct of the corresponding maps in $\Delta$
$$
(\lambda, \mu) = \lambda \lsqcup{} \mu \in \Hom_\Delta([m+n], [p+q]).
$$
Recall that a refinement $\nu$ of $(\lambda, \mu)$ is represented by a diagram $[m+n] \to [s] \epi [p+q]$. The restriction
$\nu|_\lambda$ is then defined as the pullback
$$
\xymatrix{
[m] \ar[r] \ar_{\nu|_\lambda}[d] & [m+n] \ar^\nu[d] \\
[s] \mathop{\times}\limits_{[p+q]} [p] \ar[r] \ar@{->>}[d] & [s] \ar@{->>}[d] \\
[p] \ar[r] & [p+q].
}
$$

For two partitions $\lambda \in \Lambda_p(m)$ and $\mu \in \Lambda_q(m)$ we define their {\em merge} as the ordered pushout
$$
\xymatrix@=4em{
[m] \ar^\lambda[r] \ar_\mu[d] \ar^-{\lambda \cup \mu}[dr] & [p] \ar^{\lambda_*(\mu)}[d] \\
[q] \ar_{\mu_*(\lambda)}[r] & \bullet
}
$$

We will also recall that there is a bijection between sets $\Lambda_q(p-1)$ and $\Lambda_p(q-1)$. This can be expressed
graphically by the following picture
{
\def\bt#1{\mathop{\bullet}\limits^{#1}}
\def\bb#1{\mathop{\bullet}\limits_{#1}}
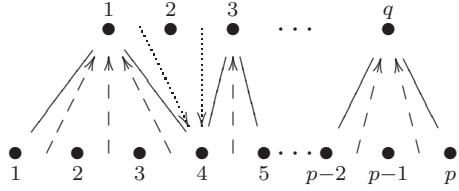
$$
\xymatrix@C=1em@R=4em@!0@H=2em{
&&& \bt{1} \ar@{-}[dlll] \ar@{-}[drrr] & *=0{} \ar@{.>}[drr] & \bt{2} & *=0{} \ar@{.>}[d] & \bt{3} \ar@{-}[dl] \ar@{-}[dr]  && \cdots & && \bt{q} \ar@{-}[dll] \ar@{-}[drr]  && \\
\bb{1} & *=0{} \ar@{-->}[urr] & \bb{2} & *=0{} \ar@{-->}[u] & \bb{3} & *=0{} \ar@{-->}[ull]& \bb{4} & *=0{} \ar@{-->}[u] & \bb{5} & \cdots & \bb{p-2} & *=0{} \ar@{-->}[ur] & \bb{p-1} & *=0{} \ar@{-->}[ul]  & \bb{p} 
}
$$
}
where dashed arrows represent an element $\mu \in \Lambda_q(p-1)$ and dotted arrows an element $\mu^* \in \Lambda_p(q-1)$. Formally, extend
$\mu$ to a map $[0, p] \cap \Z \to [0, q+1] \cap \Z$ by setting $\mu(0) = 0$ and $\mu(p) = q+1$. Define
$$
\mu^*(i) = j, \quad \text{if and only if}\quad \mu(j-1) \le i < \mu(j).
$$

Applying this construction to surjective maps we also obtain the following bijection fitting into the commutative diagram
$$
\xymatrix@R=4em@C=8em@!0{
\mathrm{Epi}_\Delta([m], [n]) \ar_{\isom}[dr] \ar^-{(-)^*}[rr] && \mathrm{Mon}_\Delta([n - 1], [m + 1]) \\
& \mathrm{Mon}_\Delta([n - 1], [m - 1]), \ar@{^(->}[ur] &
}
$$
where $[m-1]$ is embedded into $[m+1]$ as the middle part. Abusing notation we will also denote this isomorphism by $(-)^*$.

\end{nparagraph}

\begin{nparagraph}
Now let $f\from A \to B$ be a morphism of associative unital algebras, we will define spaces $D_n(f)$ of {\em differential
operators at $f$}.

Consider $\lambda_\bullet \in \Lambda_d^\circ(n)$ and $\mu_\bullet \in \Lambda_2^\circ(\lambda_i)$ for some $\lambda_i \ge 2$.
We set $\lambda' = \mu \circ_i \lambda$, the partition of $n$ of size $d + 1$. To this composition we can associate two maps, the first
map is
$$
\xymatrix@C=9em{
\Hom(A^{\tensor (d+1)}, B^{\tensor (d+1)}) \ar^-{1\tensor \cdots \tensor m_B \tensor \cdots \tensor 1}[r] & \Hom(A^{\tensor (d+1)}, B^{\tensor d}),
}
$$
where $m_B\from B \tensor B \to B$ is the multiplication in $B$, placed in $i$'th position. The second map is
$$
\xymatrix@C=9em{
\Hom(A^{\tensor d}, B^{\tensor d}) \ar^-{1\tensor \cdots \tensor \ad_m \tensor \cdots \tensor 1}[r] & \Hom(A^{\tensor (d+1)}, B^{\tensor d}),
}
$$
where $\ad_m$ sits in $i$'th position and represents the Hochschild differential (commutator with multiplication)
$\ad_m\from \Hom(A, B) \to \Hom(A \tensor A, B)$, where $A$-bimodule structure on $B$ is given via map $f$:
$$
\ad_m(g)(x, y) = g(xy) - f(x)g(y) - g(x)f(y).
$$
\end{nparagraph}

\begin{definition}
\label{def.diff.op}
The space of differential operators $D_n(f)$ for $n \ge 1$ is the kernel
$$
D_n(f) \ =\ \Ker \left( \xymatrix@C=3em{
\displaystyle \bigoplus_{\Lambda_d^\circ(n)} \Hom(A^{\tensor d}, B^{\tensor d}) \ar^-{m_B}@<0.3em>[r] \ar_-{\ad_m}@<-0.3em>[r] &
\displaystyle \bigoplus_{\Lambda_{d+1}^\circ(n) \, \downarrow\, \Lambda_d^\circ(n)} \Hom(A^{\tensor (d+1)}, B^{\tensor d})
}\right),
$$
where both sums are taken over all values of $d$. The second sum is over all arrows in the category $\ul\Lambda^\circ(n)$ from
an object in $\Lambda_{d+1}^\circ(n)$ to an object in $\Lambda_d^\circ(n)$, i.e. the set of all compositions $\lambda' = \mu \circ_i \lambda$
as above. The top and bottom arrows are the sums of all maps of the first and second type respectively described above.

\end{definition}

For an operator $P \in D_n(f)$ we will denote $P_\lambda$ the projection of $P$ to the $\lambda$-component of the direct sum.
We can extend $P$ to an element of
$$
P \in \prod_{\Lambda_d(n)} \Hom(A^{\tensor d}, B^{\tensor d})
$$
that will be also denoted by $P$,
in the following way. Let $\lambda' = (\lambda_1, \ldots \lambda_i, 0, \lambda_{i+1}, \ldots \lambda_d)$, and
$$
P_\lambda(a_1, \ldots, a_d) = \sum b^{(1)} \tensor \cdots \tensor b^{(d)},
$$
we define
$$
P_{\lambda'}(a_1, \ldots, a_i, a, a_{i+1},\ldots, a_d) = \sum b^{(1)} \tensor \cdots b^{(i)} \tensor f(a) \tensor b^{(i+1)} \tensor \cdots \tensor b^{(d)}.
$$

By convention we extend the definition of $D_n(f)$ to $n = 0$ by setting $D_0(f) = k\cdot f$, the $k$-linear span of the
map $f\from A \to B$.

\begin{proposition}
For any $P \in D_n(f)$ the commutator $[m, P] = 0$, in the sense that for every $d \ge 1$ the following holds:
$$
P_{(n)}(a_1 \cdots a_d) = \sum_{\lambda \in \Lambda_d(n)} m_B \circ P_\lambda(a_1, \ldots , a_d),
$$
where $m_B\from B^{\tensor d} \to B$ stands for the iterated product in $B$.
\end{proposition}

\proof The statement follows immediately from the definition by induction on $d$. Indeed, for $d = 1$ it is clear.
So, assuming that it is true for $d$ let us establish it for $d+1$. We have
$$
P_{(n)}(a_1 \cdots a_d a_{d+1}) = \sum_{\lambda \in \Lambda_d(n)} m_B \circ P_\lambda(a_1, \ldots , a_d a_{d+1}),
$$
and applying the definition this sum can be rewritten as
$$
\sum_{\lambda, \mu} m_B \circ P_{\mu \circ_d \lambda}(a_1, \ldots , a_d, a_{d+1}),
$$
over all $\lambda \in \Lambda_d(n)$ and $\mu \in \Lambda_2(\lambda_d)$. This in turn is clearly the same as taking the
sum over all $\lambda \in \Lambda_{d+1}(n)$.

\qed

\begin{proposition}
\label{sub-operators}
Let $P \in D_n(f)$, then for any $\lambda \in \Lambda_d(n)$ and elements $a_j \in A$, $1 \le j \le (d-1)$, the collection $\{Q_\mu\}$
for $\mu \in \Lambda_k(\lambda_i)$, defined by
$$
Q_\mu(b_1, \ldots b_k) = P_{\mu \circ_i \lambda}(a_1, \ldots, a_{i-1}, b_1, \ldots b_k, a_i, \ldots a_{d-1})
$$
is an element of $D_{\lambda_i}(f) \tensor B^{\tensor (d-1)}$.
\end{proposition}

\proof
This follows from the form of the complex defining $D_n(f)$, once one observes that the category of refinements $\lambda'$ of $\lambda$,
that leave all intervals other than $i$'th interval unchanged, is equivalent to the category of refinements $\mu$ of the
trivial partition $(\lambda_i)$. The equivalence is given by $\mu \mapsto \lambda' = \mu \circ_i \lambda$.

\qed

\begin{example}
Let $f\from A \to A$ be the identity map $\id_A$, in this case we write $D_n(A) := D_n(\id_A)$. From the proposition above
we see that the first order operators satisfy
$$
P_1(ab) = P_1(a)b + a P_1(b),
$$
i.e. they are precisely the derivations of $A$. The second order operators satisfy
$$
P_2(ab) = P_2(a) b + a P_2(b) + m_A \circ P_{1,1}(a,b),
$$
where $P_{1,1}\from A \tensor A \to A \tensor A$ is a derivation in both the first and second argument, with respect to the
$A$-bimodule structures on the first and the second factors respectively.

\end{example}

\begin{nparagraph}[Multi-differential operators.]
First we extend the notion of differential operators (\ref{def.diff.op}) to ``operators with $p$ inputs and $p$ outputs''.
Of course, since our differential operators are already collection of maps with various number of inputs and outputs this
should be understood as referring to the leading component of an operator.

For $\lambda \in \Lambda_p(n)$ we denote by $\ul\Lambda^\circ(\lambda) \subset (\ul\Lambda(n) \downarrow \lambda)$ the full
subcategory spanned by non-degenerate refinements of $\lambda$. It is easy to see that for two such refinements $\mu, \mu'$
the set $\Hom_{\ul\Lambda^\circ(\lambda)}(\mu', \mu)$ consists of a single element if $\mu'$ is a refinement of $\mu$
(necessarily non-degenerate) and empty otherwise. We write $\Lambda_d^\circ(\lambda)$ for the set of objects $\mu \in \ul\Lambda^\circ(\lambda)$,
such that $\mu$ is a partition of size $d$.

\begin{definition}
\label{def.multidiff}
The space of differential operators $D_\lambda(f)$ of order $\lambda \in \Lambda_p(n)$ is the kernel
$$
D_\lambda(f) \ =\ \Ker \left( \xymatrix@C=3em{
\displaystyle \bigoplus_{\Lambda_d^\circ(\lambda)} \Hom(A^{\tensor d}, B^{\tensor d}) \ar^-{m_B}@<0.3em>[r] \ar_-{\ad_m}@<-0.3em>[r] &
\displaystyle \bigoplus_{\Lambda_{d+1}^\circ(\lambda) \, \downarrow\, \Lambda_d^\circ(\lambda)} \Hom(A^{\tensor (d+1)}, B^{\tensor d})
}\right).
$$
\end{definition}

For example, if $P \in D_n(f)$, then for any $\lambda \in \Lambda(n)$ the collection $\{P_\mu \mid \mu \in \ul\Lambda^\circ(\lambda)\}$
is an element of $D_\lambda(f)$.
\end{nparagraph}

\begin{nparagraph}
\label{multi.qp}
Next we extend this to ``operators with $q$ inputs and $p$ outputs'', with $p \le q$. Informally, such an operator is an
operator with $q$ inputs and outputs, composed with multiplication $m_\pi\from B^{\tensor q} \to B^{\tensor p}$ according
to some partition $\pi$ of $q$ into $p$ parts.

Consider refinement $\lambda' \in \Lambda_{q'}(n)$ of $\lambda \in \Lambda_q(n)$ and let $\rho\from [q'] \to [q]$ be the
corresponding map. For any partition $\pi \in \Lambda^\circ_p(q)$ we form the lift $\pi' = \pi'(\rho) \in \Lambda^\circ_{p+q'-q}(q')$ as follows.
We have two injective maps
$$
\xymatrix{
[p-1] \ar@{^(->}^{\pi^*}[r] & [q-1] \ar@{^(->}^{\rho^*}[r] & [q'-1],
}
$$
and let $i$ be the inclusion of the union
$$
\xymatrix{
\Im(\pi\rho)^* \cup ([q'-1] - \Im \rho^*) \ar@{^(->}^-{i}[r] & [q'-1].
}
$$
Since the size of the union is $(p - 1 + q' - q)$ applying the $*$-bijection again we obtain
$$
\pi' = i^* \from [q'] \epi [p + q' - q].
$$
In other words $\pi'$ can be defined as the unique partition in $\Lambda^\circ_{p+q'-q}(q')$ that satisfies the two conditions
\begin{itemize}
\item for every $x \in [q]$ all elements in $\rho^{-1}(x)$ belong to different parts of partition $\pi'$,
\item $\pi$ is the merge of $\pi'$ along $\rho$, i.e. the diagram
$$
\xymatrix{
[q'] \ar^\rho[r] \ar_{\pi'}[d] & [q] \ar_{\pi}[d] \\
[p+q'-q] \ar[r] & [p]
}
$$
is an ordered pushout.
\end{itemize}
Let us illustrate this definition by an example, let $\pi = (3)$ be the trivial partition and $\rho\from [6] \to [3]$,
$\rho = (1, 3, 2)$. Then $\pi' = (2, 1, 2, 1) \from [6] \to [4]$ as shown in the picture below. The ligatures on the
right hand side represent partition $\rho$.
\vskip 1em
{
\def\b{*=0 {\bullet}}
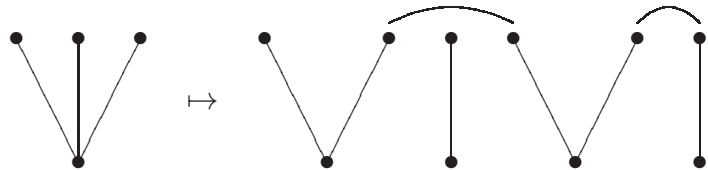
$$
\xymatrix@C=2em@R=2em@!0@H=2em{
\b \ar@{-}[ddr] & \b \ar@{-}[dd] & \b \ar@{-}[ddl] & & \b \ar@{-}[ddr] && \b \ar@{-}[ddl] \ar@{-}@<0.5em>@/^/[rr]& \b \ar@{-}[dd]& \b \ar@{-}[ddr] && \b \ar@{-}[ddl] \ar@{-}@<0.5em>@/^/[r] & \b \ar@{-}[dd] \\
&  &  & \mapsto & & && & && &\\
& \b & & &  & \b & & \b & & \b & & \b
}
$$
}
We also observe, that if $\pi = 1^q = (1, \ldots, 1)$ the finest non-degenerate partition of $q$, then $\pi'(\rho) = 1^{q'}$
for any $\rho$.

For every $\pi \in \Lambda^\circ_p(q)$ we write $\Hom(A^{\tensor d}, B^{\tensor q})_\pi$ for the copy of the $\Hom$-space indexed by $\pi$.
The map $\ad_m$ goes from $\pi$-component to $\pi$-component, and for every $\rho\from [q'] \to [q]$ the map $m_B$ goes
$$
m_B\from \Hom(A^{\tensor q'}, B^{\tensor q'})_{\pi'(\rho)} \to \Hom(A^{\tensor q'}, B^{\tensor q})_{\pi}.
$$

\end{nparagraph}

\begin{definition}
The space of differential operators $D_{\lambda,\pi}(f)$ of order $\lambda \in \Lambda_q(n)$ and output type $\pi \in \Lambda_p(q)$
is the kernel
$$
D_{\lambda,\pi}(f) \ =\ \Ker \left( \xymatrix@C=3em{
\displaystyle \bigoplus_{\Lambda_d^\circ(\lambda)} \Hom(A^{\tensor d}, B^{\tensor d})_{\pi'(\rho)} \ar^-{m_B}@<0.3em>[r] \ar_-{\ad_m}@<-0.3em>[r] &
\displaystyle \bigoplus_{\Lambda_{d+1}^\circ(\lambda) \, \downarrow\, \Lambda_d^\circ(\lambda)} \Hom(A^{\tensor (d+1)}, B^{\tensor d})_{\pi'(\rho)}
}\right),
$$
where $\rho\from [d] \to [q]$ is the map corresponding to a refinement $\lambda' \in \Lambda^\circ_d(\lambda)$.
\end{definition}

Of course this complex is the same as in definition of $D_\lambda$ for any $\pi$, and the decoration of components
are purely for bookkeeping purposes that will be later used in definition of the composition.

\begin{nparagraph}
Finally, we extend the notion of differential operators to ``$q$ inputs and $p$ outputs'' for any $p, q \ge 1$.
Let $A^e = A \tensor A^\op$ be the enveloping algebra of $A$ and set $B = T^\bullet_{A} (A^e)$. It can be written as
the direct sum
$$
B = A \oplus A^{\tensor 2} \oplus A^{\tensor 3} \oplus \ldots,
$$
and let $f\from A \to B$ be the unit map, i.e. the embedding of $A$ as the first summand.
We will write
$$
\cD_n(A) = \bigoplus_{\lambda \in \ul\Lambda(n) \atop \pi \in \ul\Lambda^\circ(|\lambda|)} D_{\lambda,\pi}(f),
$$
and call it the space of {\em multi-differential operators} on $A$.

\end{nparagraph}

\begin{nparagraph}[Bigrading on $\cD_n(A)$.]
First let us introduce a grading on $B$ by putting $B_n = A^{\tensor n+1}$. The multiplication $m_B$ is of degree $0$
with respect to this grading:
$$
m_B\from A^{\tensor m+1} \tensor A^{\tensor n+1} \to A^{\tensor (m+n+1)}.
$$
We will also equip $B^{\tensor d}$ with the total grading
$$
(B^{\tensor d})_n \ =\  \bigoplus_{n_1 + \cdots + n_d = n} B_{n_1} \tensor \cdots \tensor B_{n_d}.
$$

Furthermore, we define a bigrading on the spaces $\Hom(A^{\tensor q}, B^{\tensor d})$ by putting $\phi \in \Hom(A^{\tensor q}, (B^{\tensor d})_p)$
in bidegree $(p, q)$. Then the map $m_B$ from the definition of $D_\lambda(f)$ is of bidegree $(0, 0)$ and the map $\ad_m$ is
of bidegree $(0, 1)$. Therefore, by restricting to bidegrees $(p, \bullet)$ we obtain a subcomplex of the two-term complex
in the definition of $D_\lambda(f)$ (\ref{def.multidiff}).
Denoting by $D_\lambda(f)_p$ the kernel of this subcomplex, we have $D_\lambda(f) = \bigoplus D_\lambda(f)_p$. In fact we
obtain even finer decomposition $D_\lambda(f)_p = \bigoplus_{\mu \in \Lambda_q(p)} D_\lambda^\mu$.

Finally, we define bigrading on $\cD_n(A)$ by putting for $p, q \ge 1$
$$
\cD_n(A)(p, q) = \bigoplus_{\lambda \in \Lambda_q(n)  \atop \pi \in \ul\Lambda^\circ(|\lambda|) } D_\lambda(f)_{p-|\pi|}.
$$
One can think of $q$ as the number of inputs of a differential operator and $p$ as the number of outputs.

\end{nparagraph}

\vskip 5em 

\subsection{Planar prop structure}
\label{sect-planar-prop}

We begin by constructing an associative algebra structure on the space $D(A) = \bigoplus D_n(A)$. Here, as before $D_n(A) = D_n(\id_A)$
and $D_0(A) = k$.

For two partitions $\lambda' \in \Lambda_d(m)$ and $\lambda'' \in \Lambda_d(n)$ we write
$$
\lambda = \lambda' + \lambda'' = (\lambda'_1 + \lambda''_1, \ldots \lambda'_d + \lambda''_d) \in \Lambda_d(m + n).
$$

Now, for two operators $P \in D_n(A)$, $Q \in D_m(A)$ we form their composition $Q \circ P$ by putting
\begin{equation}
\label{comp.DA}
(Q \circ P)_\lambda = \sum_{\lambda = \lambda' + \lambda''} Q_{\lambda'} \circ P_{\lambda''},
\end{equation}
with $\lambda' \in \Lambda_d(m)$ and $\lambda'' \in \Lambda_d(n)$.

\begin{proposition}
\label{assoc.DA}
The space  $(D(A), \circ)$ is an associative algebra.
\end{proposition}

\proof
It is clear from the definition that $\circ$ is an associative operation and that $1 \in D_0(A)$ is the unit. It remains to
show that it is a well defined operation on $D(A)$. We have
$$
P_{\lambda''}(\ldots, a_i a_i', \ldots) = \sum_{\mu'' \in \Lambda_2(\lambda''_i)} m_B^{(i)} \circ P_{\mu'' \circ_i \lambda''}(\ldots, a_i, a_i', \ldots),
$$
where $m_B^{(i)}$ is the product of elements in positions $i$ and $(i+1)$. Furthermore,
$$
Q_{\lambda'} \circ P_{\lambda''}(\ldots, a_i a_i', \ldots) =
\sum_{\mu' \in\Lambda_2(\lambda'_i)} \sum_{\mu'' \in \Lambda_2(\lambda''_i)} m_B^{(i)} \circ Q_{\mu' \circ_i \lambda'} \circ P_{\mu'' \circ_i \lambda''}(\ldots, a_i, a_i', \ldots).
$$
It remains to observe that the set of pairs $(\nu', \nu'')$, such that $\nu' + \nu'' = \mu \circ_i \lambda$ with $\mu$ ranging over
all elements in $\Lambda_2(\lambda_i)$ is in bijection with the set of tuples $(\lambda', \lambda'', \mu', \mu'')$ as above.
This shows that $(Q \circ P) \in D_{m+n}(A)$.

\qed

The rest of this section is dedicated to proving the following theorem.

\begin{theorem}
\label{thm.planar.prop}
The space of multi-differential operators $\cD(A) = \bigoplus \cD_n(A)$ has a structure of a reduced planar prop.
\end{theorem}

\proof
First we restrict our consideration to the output types $\pi = (1, \ldots 1)$.
According to proposition (\ref{elem.prop}) it is enough to construct operations of horizontal and vertical composition
and show that they satisfy relations in definition (\ref{def.elem.prop}). We put the horizontal unit $u^0 = 1 \in k = \cD(0, 0)$
and the vertical unit $u = \id_A \in \cD_0(1, 1)$.

For $P \in D_\lambda^{\lambda'}$ and $Q \in D_\mu^{\mu'}$ we define the horizontal composition in the following way. For any non-degenerate
refinement $\nu$ of $(\lambda, \mu)$ we put
$$
(P \circ_h Q)_\nu = P_{\nu|_\lambda} \tensor Q_{\nu|_\mu}.
$$
It is clear that this defines an element in $D_{(\lambda, \mu)}^{(\lambda', \mu')}$. In other words we obtain an operation
$$
\circ_h \from \cD_n(p, q) \tensor \cD_{n'}(p', q') \to \cD_{n+n'}(p+p', q+q').
$$

Next we define the vertical composition. Let $P \in D_\lambda^\mu$ and $Q \in D_\nu$, for some $\lambda \in \Lambda_q(n)$,
$\mu \in \Lambda_q(p-1)$ and $\nu \in \Lambda_{p+q-1}(m)$. Recall that we can think of $\lambda$ and $\nu$ as partitions of the
order of operators $P$ and $Q$ among their inputs and $\mu$ as the partition of outputs of $P$. Consider the following
diagram.

\begin{equation}
\label{vcomp.diag}
\xymatrix@C=4em{
& [n] \ar_-{\lambda}[dl] \ar_-{\lambda'}[d] \ar@{^(->}[r] & [m+n] \ar[d] & [m] \ar^-{\lambda''}[d] \ar@{_(->}[l] \ar^{\nu}[dr]  &  \\
[q] & [q'] \ar@{->>}^-{\rho}[l] \ar@{^(->}^-{i_1}[r] \ar_-{=}[dr] & [q'] \lsqcup{[q']} [q' + p - 1] \ar[d] & [q' + p - 1] \ar^-{\id_{[q']} \lsqcup{} \mu'}[dl] \ar@{_(->}_-{i_2}[l] \ar@{->>}_-{\alpha}[r] & [q + p - 1] \\
& & [q'] & & 
}
\end{equation}
Here $\mu' \in \Lambda_{q'}(p-1)$ is a refinement of $\mu$, such that $\mu = \rho \circ \mu'$, maps
$i_1$ and $i_2$ are the structure maps of the ordered coproduct over $[q']$, and the map $\alpha$ is defined as follows.
As was discussed before a surjective map $\rho\from [q'] \to [q]$ determines an injective map $\rho^*\from [q-1] \to [q'-1]$.
Applying the $*$-duality to $\mu'\from [p-1] \to [q']$ fiberwise over $[q]$ we construct $(\mu'/\mu)^*$:
$$
(\mu'/\mu)^* \from ([q'-1] - \Im\rho^*) \to [q + p - 1].
$$
Then we have
$$
[q' + p - 1] \isom [q+p-1] \lsqcup{[q+p-1]} \Im (\mu'/\mu)^*,
$$
and $\alpha$ is the structure map of the ordered coproduct.

Denote the composition in the middle column by $\tau$. We define
for any $\tau$
$$
(Q \circ_v P)_\tau = \sum_{\lambda', \lambda''} Q_{\lambda''} \circ P_{\lambda'},
$$
where the sum is taken over all $\lambda', \lambda''$ that fit into a diagram (\ref{vcomp.diag}) (if there are no such diagrams
then the $\tau$-component of the composition is $0$).

Notice that the composition formula in (\ref{comp.DA}) is a special case of the vertical composition when $p = q = 1$. The proof
that this composition is a well-defined operator as well as the associativity of the composition is analogous to the proof of
proposition (\ref{assoc.DA}), and is done by establishing bijection between two sets of diagrams.

This gives us an operation
$$
\circ_v \from \cD_m(r, p) \tensor \cD_n(p, q) \to \cD_{m+n}(r, q).
$$

It is clear that $u^0$ is the unit for the horizontal composition. To show that the compatibility condition of definition (\ref{def.elem.prop})
holds it is enough to show that powers of $u$ are units for the vertical composition. Indeed, let $P = u^q$, $Q \in \cD_m(p, q)$,
the component $P_{\lambda'}$ is non-zero if and only if the refinement $\lambda' = (0, \ldots, 0)$. Furthermore, we have
$p = 1$, and $\alpha = \rho$. Therefore the only non-zero contribution to the sum comes from $\tau = \lambda''$ and since
$\lambda''$ ranges over all refinements of $\nu$ we have $(Q \circ_v u^q) = Q$.

Similarly, for $P \in \cD_n(p, q)$, $Q = u^p$, we have $\lambda'' = (0, \ldots, 0)$, and therefore $\tau = \lambda'$ ranges
over all refinements of $\lambda$, so we have $(u^p \circ_v P) = P$.

\begin{nparagraph}[Other output types.]
\label{output.types}
For $P \in D_{\lambda,\pi}$ and $Q \in D_{\nu, \sigma}$ we define their horizontal composition $P \bullet_h Q$ by taking the horizontal
composition of untyped operators $P \circ_h Q$ as constructed before and putting it in the output type $(\pi, \sigma)$:
$$
P \bullet_h Q = (P \circ_h Q)[(\pi, \sigma)].
$$
The associativity of this operation follows from the associativity of operations $\circ_h$ above and concatenation of partitions.

Recall that for any refinement $\nu'$ of $\nu$ the collection $Q_{(\nu')} = \{Q_{\nu''} \mid \nu'' \in \ul\Lambda^\circ(\nu')\}$ is
an element of $D_{\nu', \sigma'(\rho)}$, where $\rho$ is the map determined by the refinement $\nu' \to \nu$.
Assume $P \in D_{\lambda,\pi}^\mu$, $\mu \in \Lambda_q(p-1)$ and define the vertical composition as
$$
(Q \bullet_v P) = \sum_{\nu=\wtilde\pi \nu'} (Q_{(\nu')} \circ_v P) [\wtilde\mu_*(\sigma) \pi],
$$
where $\wtilde\pi = \pi'(\id_{[q]} \sqcup \mu)$ (see paragraph \ref{multi.qp} for the construction), $\wtilde\mu_*(\sigma)$ is the
ordered pushout of $\sigma$ along $\wtilde\mu$ and $\wtilde\mu \from [p+s-1] \to [s]$
with $s = |\pi|$ is defined as the structure map of the ordered coproduct
$$
\xymatrix@C=6em{
[s] \lsqcup{[s]} [p-1] \ar^-{\id_{[s]} \sqcup (\pi \mu)}[r] & [s]
}
$$

Let us show associativity of this operation. This follows from a combination of associativity of $\circ_v$, the fact that
$(Q_{(\nu')})_{(\nu'')} = Q_{(\nu'')}$ (compatibility of restrictions) and associativity of the output type decoration,
that we will now check.

Let $P \in D_{\lambda,\pi_P}^{\mu_P}$, $Q \in D_{\lambda,\pi_Q}^{\mu_Q}$ and $R \in D_{\lambda,\pi_R}^{\mu_R}$,
we write $s_P = |\pi_P|$, $s_Q = |\pi_Q|$ and $s_R = |\pi_R|$, also let $\mu_P \in \Lambda_q(p_P-1)$ and $\mu_Q \in \Lambda_{p + s_P - 1}(p_Q - 1)$.
The associativity of decoration follows from commutativity
of the diagram below, since the decoration of the composite $R \bullet_v Q \bullet_v P$ is represented by the top row,
and all three squares are ordered pushouts.

$$
\xymatrix{
[q] \ar^-{\pi_P}[r] & [s_P] \ar[r] & |\wtilde\mu_P \cup \pi_Q| \ar[r] & \bullet \\
& [p_P + s_P - 1] \ar^-{\wtilde\mu_P}[u] \ar^-{\pi_Q}[r] & [s_Q] \ar[u] \ar[r] & |\wtilde\mu_Q \cup \pi_R| \ar[u] \\
& & [p_Q + s_Q - 1] \ar^-{\wtilde\mu_Q}[u] \ar^-{\pi_R}[r] & [s_R]. \ar[u]
}
$$
\vskip 1em

\end{nparagraph}

This completes the proof of theorem (\ref{thm.planar.prop}).

\qed

\begin{nparagraph}[Modular prop structure.]
In addition to gradings of $\cD(A)$ by the number of inputs and outputs we can define another grading that we will call the
genus grading. The {\em genus} of an operator $P \in \cD_n(p, q)$ is defined as
$$
g(P) = n - p + 1.
$$
In general $g(P)$ can be any integer.
\end{nparagraph}

\begin{theorem}
The space of multi-differential operators $\cD(A)$ equipped with the genus grading has a structure of a modular planar prop.
\end{theorem}

\proof
We need to check that this notion of genus is compatible with the composition along planar graphs as explained in
(\ref{genus.markings}). In fact it is enough to check it for horizontal and vertical compositions. Let $P \in \cD_n(p, q)$
and $Q \in \cD_m(r, s)$, it is clear that for the horizontal composition
$$
g(P \bullet_h Q) = m + n - p - r + 1 = g(P) + g(Q) - 1,
$$
and the genus of the graph with two disjoint contractible components is $(-1)$. For the vertical composition we have
$$
g(Q \bullet_v P) = m + n - r + 1 = g(P) + g(Q) + (p - 1),
$$
and the genus of the graph for this vertical composition is $(p - 1)$.

\qed

\begin{definition}
We say that an element $P \in \cD(A)$ is of totally positive genus if for every non-zero component $P_\lambda^\mu \in D_\lambda^\mu$
we have $\lambda \ge \mu$ (termwise). We will denote the subprop of operators of totally positive genus by $\cD^{\ge 0}(A)$.
\end{definition}

\vskip 5em
\subsection{Automorphisms of associative families}
\label{sect-assoc-families}

In this section we construct another planar prop $\cE(A)$ for an associative algebra $A$ and a morphism of props $\cE(A) \to \cD(A)$.

Let $R = k\<h_1, \ldots, h_n\>$ be the free associative algebra generated by the set $H = \{h_1, \ldots, h_n\}$,
and consider the coproduct in the category of associative algebras $A * R$. Denote by $H^*$ the set of words in alphabet
$H$, and for any word $w \in H^*$ we write $|w|$ for its length.
Elements of $A*R$ have the form $a_1 h_{i_1} a_2 h_{i_2} \ldots h_{i_k} a_{k+1}$, for some $a_j \in A$.
In other words
$$
A*R = \bigoplus_{w \in H^*} A^{\tensor |w|+1}.
$$
Equip $A*R$ with the decreasing filtration by the length of $w$
$$
F_i(A*R) = \bigoplus_{|w| \ge i} A^{\tensor |w|+1}
$$
and denote $A \sh R$ the completion of $A*R$ with respect to this filtration.

One can think of $A \sh R$ as the trivial formal associative family of algebras over $R$ with fiber $A$. Consider the group
$\Aut_R(A \sh R)$ of automorphisms of $A \sh R$ as an algebra over $R$. By restricting any such automorphism $\phi$ to the ``fiber
over $0 \in \Spec R$'' we obtain an automorphism of $A$:
$$
\xymatrix@C=4em{
A \ar@{^(->}[r] & A \sh R \ar^\phi[r] & A \sh R \ar^-{h_i \mapsto 0}[r]  & A.
}
$$
We consider the subgroup $\Aut^\circ_R(A \sh R)$ of automorphisms such that their restriction to $A$ is identity. Explicitly
any such automorphism is given by a collection of maps
$$
\phi = \{\phi_w \from A \to A^{\tensor |w|+1} \mid w \in H^*\},
$$
such that $\phi_{\epsilon}(a) = a$, for the empty word $\epsilon$ and satisfying
\begin{equation}
\label{phi.aut}
\phi_w(ab) = \sum_{w = w'w''} \phi_{w'}(a) \phi_{w''}(b).
\end{equation}

Fix $w \in H^*$ and let $n = |w|$. For any partition $\lambda \in \ul\Lambda^\circ(n)$ we set $w_\lambda(i)$ to be the product of letters
of $w$ in the $i$'th part of the partition.

\begin{lemma}
For any $w \in H^*$ the collection $\{\phi_\lambda \mid \lambda \in \ul\Lambda^\circ(n)\}$ defined as
$$
\phi_\lambda = \phi_{w_\lambda(1)} \tensor \cdots \tensor \phi_{w_\lambda(|\lambda|)}
$$
is an element of $\cD_n(A)$.
\end{lemma}

\proof
Clearly, the fact that $\left(\sum \phi_\lambda\right) \in \Ker(\ad_m - m_B)$ immediately follows from conditions (\ref{phi.aut}).

\qed

\begin{definition}
Denote by $\cE(A)$ the planar prop generated by elements
$$
\{\phi_w \in \cE(|w|+1, 1) \mid \phi \in \Aut^\circ_R(A \sh R), \ w \in H^*\}
$$
and $m_A \in \cE(1, 2)$ modulo relations (\ref{phi.aut}) and associativity of $m_A$.
\end{definition}

\begin{proposition}
There is a morphism of planar props $r\from \cE(A) \to \cD(A)$ that sends each $\phi_w$ to the collection $\{\phi_\lambda\}$ as
in the lemma and $m_A$ to $(\id_A \circ_h \id_A)[(2)]$ (for the notation see paragraph \ref{output.types}).
\end{proposition}

\proof
We just need to check that $\{\phi_\lambda\}$ and $(\id_A \circ_h \id_A)[(2)]$ satisfy relations (\ref{phi.aut}) in $\cD(A)$.
But this immediately follows from the definition of the vertical composition:
$$
\{\phi_\lambda\} \bullet_v (\id_A \circ_h \id_A)[(2)]\ =\ \sum_{\lambda = (\lambda', \lambda'')}
(\{\phi_{\lambda'}\} \circ_h \{\phi_{\lambda''}\}) \circ_v (\id_A \circ_h \id_A) [(2)].
$$

\qed

We will also need the following two simple constructions in the proof of theorem (\ref{thm.surjective}).

\begin{lemma}
\label{deriv.aut}
Let $H$ be a finite set and $R = k\<H\>$. Any collection of derivations
$$
\{\d_h \from A \to A \tensor A \mid h \in H\}
$$
determines an automorphism $\phi \in \Aut^\circ_R(A \sh R)$ by putting for every $w = h_{i_1} h_{i_2} \ldots h_{i_k} \in H^*$
$$
\xymatrix@C=4.5em{
\phi_w \from A \ar^{\d_{i_1}}[r] & A \tensor A \ar^-{1\tensor \d_{i_2}}[r] & A \tensor (A \tensor A) \ar^-{1 \tensor 1 \tensor \d_{i_3}}[r] & \ldots \ar^-{1^{\tensor (k-1)} \tensor \d_{i_k}}[r] & A^{\tensor |w|+1}.
}
$$
\end{lemma}

\proof
The relations (\ref{phi.aut}) are immediately verified by using the fact that each $\d_h$ is a derivation.

\qed

\begin{lemma}
\label{degen.aut}
Let $\sigma \from G \epi H$ be a surjective map in the category $\Delta$. It determines a map $s = s(\sigma) \from \Aut^\circ_{k\<H\>}(A \sh k\<H\>) \to \Aut^\circ_{k\<G\>}(A \sh k\<G\>)$.
\end{lemma}

\proof
Map $\sigma$ determines an algebra map $\sigma^* \from k\<H\> \to k\<G\>$ by sending each generator $h \in H$ to
$$
\sigma^*(h) = \prod_{g \in \sigma^{-1}(h)} g,
$$
where the product is taken in the order prescribed by the total order in $G$. The diagram
$$
\xymatrix@=4em{
A \sh k\<H\> \ar^\phi[r] \ar_{\sigma^*}[d] & A \sh k\<H\> \ar^{\sigma^*}[d] \\
A \sh k\<G\> \ar@{-->}_{s(\phi)}[r] & A \sh k\<G\>
}
$$
can be completed by the dashed arrow in a unique way, since it has to be an automorphism over $k\<G\>$.

\qed

\vskip 5em
\subsection{Filtration and the symbol map}
\label{sect-symbol}

Denote by $\Delta^{>0}_{\mathrm{Epi}}$ the subcategory of $\Delta$ formed by objects $\{[n]\}$ with $n > 0$ and surjective
maps between them. For any category $\C$ an {\em epi-simplicial} object in $\C$ is a functor $(\Delta^{>0}_{\mathrm{Epi}})^\op \to \C$,
in other words it is the part of a simplicial object formed by all degeneracy maps.

Consider spaces $D_n(f)$ for some morphism of associative algebras $f\from A \to B$. To any surjective map $\sigma\from [m] \epi [n]$,
we associate a degeneracy map $s = s(\sigma)\from D_n(f) \to D_m(f)$, by putting for any $P \in D_n(f)$
$$
s(P)_{\lambda'} = \begin{cases}
P_\lambda,&\text{if $\lambda' = \lambda \sigma$}\\
0,& \text{otherwise}.
\end{cases}
$$
It is straightforward to check that this collection is a well defined element $s(P) \in D_m(f)$. Applying this to $B = T^\bullet_A (A^e)$
we obtain

\begin{proposition}
The collection $\{\cD_{n-1}(A) \mid n \ge 1\}$ and degeneracy maps $s$ form an epi-simplicial vector space.
\end{proposition}

\qed

Using degeneracy maps we construct the order filtration on $\cD_m(A)$ by putting
$$
F_n(\cD_m(A)) = \sum_{\sigma\from [m] \epi [n]} \Im \left(s(\sigma)\from \cD_n(A) \to \cD_m(A) \right).
$$

\begin{nparagraph}
Recall that an associative algebra $A$ is called {\em formally smooth} if the space of non-commutative forms $\Omega^1_A$,
defined as the kernel
$$
\xymatrix{
\Omega^1_A \ar@{^(->}[r] & A \tensor A \ar@{->>}^-{m}[r] & A,
}
$$
is a finitely generated projective $A^e$-module. Recall also that for a formally smooth algebra $A$ the Hochschild cohomology
groups $HH^i(A, M) = 0$, for all $i \ge 2$, and any $A$-bimodule $M$.
\end{nparagraph}

\begin{lemma}
\label{lemma.2diff}
Let $A$ be a formally smooth algebra and $f\from A \to B$ any morphism of algebras. Then any $P \in D_{(1,1)}(f)$ can be lifted
to an element $Q \in D_2(f)$, such that $Q_{(1,1)} = P$.
\end{lemma}

\proof
Consider the following diagram
$$
\xymatrix@R=3em{
\Hom(A, B) \ar^-{d}[r] & \Hom(A^{\tensor 2}, B) \ar^-{d}[r] & \Hom(A^{\tensor 3}, B) \ar^-{d}[r] & \cdots \\
& \Hom(A^{\tensor 2}, B^{\tensor 2}) \ar@{-->}[ul] \ar_-{m_B}[u] & &
}
$$
The top row is the Hochschild complex, where $A$-bimodule structure on $B$ is given by the map $f$. Let us show that
the biderivation $P \in \Hom(A^{\tensor 2}, B^{\tensor 2})$ is sent to a Hochschild cocycle $\bar P$ by $m_B$. Indeed,
\begin{multline*}
a \bar P(b, c) - \bar P(ab, c) + \bar P(a, bc) - \bar P(a, b) c = \\
a \bar P(b, c) - (a \bar P(b, c) + \bar P_{(1,0,1)}(a, b, c))  + (\bar P_{(1,0,1)}(a, b, c) + \bar P(a, bc)) - \bar P(a, b) c = 0.
\end{multline*}
Since $HH^2(A, B) = 0$ we find that there exists an element $Q \in \Hom(A, B)$, such that $dQ = \bar P$. By construction the collection
$\{Q, P\}$ is an element of $D_2(f)$.

\qed

\begin{lemma}
\label{lemma.cocycle}
Let $A$ be any associative algebra and $P \in D_n(f)$, then the sum
$$
Z = \sum_{\Lambda^\circ_2(n)} \bar P_{\lambda} \in \Hom(A^{\tensor 2}, B)
$$
is a Hochschild $2$-cocycle.
\end{lemma}

\proof
As in the proof of the previous lemma one checks that
$$
a Z(b, c) - Z(ab, c) + Z(a, bc) - Z(a, b) c = -\sum_{\Lambda^\circ_3(n)}  \bar P_{\lambda}(a, b, c) + \sum_{\Lambda^\circ_3(n)}  \bar P_{\lambda}(a, b, c) = 0.
$$

\qed

In fact, the same argument shows that for any even $p$ (but not odd!) the sum
$$
Z_p = \sum_{\Lambda^\circ_p(n)} \bar P_{\lambda} \in \Hom(A^{\tensor p}, B)
$$
is a Hochschild $p$-cocycle.

\begin{lemma}
\label{lemma.merge}
Let $A$ be a formally smooth algebra and $P \in D_{(n,m)}(f)$. Then there exists $Q \in D_{n+m}(f)$ such that $Q_{((n,m))} = P$.
\end{lemma}

\proof
Take $\mu \in \Lambda^\circ(m+n)$. If $\mu$ is a refinement of the concatenation $((n), (m))$ then $P_\mu$ is well defined
and we set $Q_\mu = P_\mu$.
Now, if $\mu$ is not a refinement, then let $\mu_j$ be the interval spanning the gap between $(n)$ and $(m)$.
We prove the statement using two-step induction: first on the size of $\mu_j$ and then on $M = \max_{i \neq j} (\mu_i)$.
The case $\mu_j = 2$ and $M = 1$ follows from lemma (\ref{lemma.2diff}). Let $\mu$ be such that $\mu_j = 2$ and $M > 1$,
and denote $d = |\mu|$. Form the poly-Hochschild complex with terms
$$
C_{(d)}^p(A, B) = \bigoplus_{\rho\from[p] \epi [d]} \Hom(A^{\tensor p}, B^{\tensor d}).
$$
By induction we have already constructed $Q_{\mu'}$ for every refinement $\mu' \to \mu$. Using lemma (\ref{lemma.cocycle})
we see that the sum
$$
\sum_{\Lambda^\circ_{d+1}(\mu)} \bar Q_{\mu'}
$$
is a cocycle in $C_{(d)}^{d+1}$. The vanishing $HH_{(d)}^{d+1}(A, B) = 0$ implies the existence of a lift in $C_{(d)}^d$
and we put $Q_\mu$ to be this lift.

Next, assume we have constructed $Q_\mu$ for all $\mu$ with $\mu_j < k$. Let $\mu$ be such that $\mu_j = k$ and $M = 1$,
we construct $Q_\mu$ applying the same argument to the refinements of $\mu$.

It is clear that the collection $\{Q_\mu\}$ obtained this way is an element of $D_{m+n}(f)$ and $Q_{((n,m))} = P$.

\qed

\begin{proposition}
\label{symbol.map}
Let $A$ be a formally smooth algebra, then for any $\lambda \in \Lambda_q(n)$ we have the short exact sequence
$$
\xymatrix@C=4em{
F_{n-1} D_\lambda(f) \ar@{^(->}[r] & D_\lambda(f) \ar@{->>}^-{\Sigma}[r] & D_{(1^n)}(f).
}
$$
\end{proposition}

\proof
The surjectivity of $\Sigma$ immediately follows by repeated application of lemma (\ref{lemma.merge}) to an element in $D_{(1^n)}(f)$.
Let us show exactness in the middle term. Without loss of generality we may assume that $\lambda = (n)$.

Let $P \in D_n(f)$ be
such that $P_{(1^n)} = 0$. This implies that for every $\mu \in \Lambda^\circ_{n-1}(n)$ the element $P_\mu$ in fact belongs to
$D_{(1^{n-1})}$. As mentioned before it can be lifted to some $Q^{(\mu)} \in D_{n-1}(f)$. Consider element
$$
R = P - \sum_{\mu \in \Lambda^\circ_{n-1}(n)} s(\mu)(Q^{(\mu)}).
$$
Clearly, we have $R_\mu$ = 0, for all $\mu$ of size $|\mu| \ge (n-1)$. Iterating this process we find that $P \in F_{n-1} D_n$.

\qed

We will refer to $\Sigma\from D_\lambda(f) \to D_{(1^n)}(f)$ as the {\em symbol} map.

\vskip 2em
\begin{theorem}
\label{thm.surjective}
For a formally smooth algebra $A$ the map $r \from \E(A) \epi \cD^{\ge 0}(A)$ is surjective.
\end{theorem}

\proof
It is enough to restrict our attention to elements in $D_n$, since any output type decoration $\pi$ can be obtained by
composing with a suitable multiplication in $\E(A)$, and any $D_\lambda$ can be obtain from several $D_n$ by taking horizontal
composition.

Let $P \in D_n$ be an operator of totally positive genus, and assume that $\Sigma(P) \neq 0$. We also observe that
since $\Omega^1_A$ is a finitely generated $A^e$-module the natural inclusion
$$
\xymatrix@C=4em{
\Der(A, B)^{\tensor n} \ar@{^(->}^-{\isom}[r] & D_{(1^n)}(f)
}
$$
is an isomorphism. So $\Sigma(P)$ can be written as a sum of tensor products of derivations, and total positivity implies
that every derivation has either one or two outputs.

If a derivation has only one output we replace it with the composition of its lift to a derivation with two outputs and multiplication.
Existence of such a lift follows from the formal smoothness of $A$ and the following exact sequence
$$
\xymatrix{
\Hom(\Omega^1_A, A \tensor A) \ar^-{m}[r] & \Hom(\Omega^1_A, A) \ar[r] & \Ext^1(\Omega^1_A, \Omega^1_A),
}
$$
since the $\Ext$ vanishes.

Let $H = [n]$ as an ordered set, and $R = k\<H\>$. Using lemma (\ref{deriv.aut}) we construct an automorphism $\phi \in \Aut^\circ_R(A \sh R)$.
Take $w = h_1 h_2 \ldots h_n$, and put $Q = r(\phi_w)$. By construction $\Sigma(Q) = \Sigma(P)$, so we reduced the question to a degeneracy
of a lower order operator.

It is clear from the definition of degeneracy maps $s(\sigma)$ and the construction of lemma (\ref{degen.aut}) that the following
square commutes
$$
\xymatrix@=4em{
\phi_w \ar@{|->}[r] \ar@{|->}_-{r}[d] & s(\sigma)(\phi)_{\sigma^*(w)} \ar@{|->}^-{r}[d] \\
Q \ar@{|->}[r] & s(\sigma)(Q).
}
$$
Therefore the image of $r$ is closed under taking degeneracy maps. We conclude the proof by using induction on the order $n$.

\qed

\vfill\eject

\end{document}